\documentclass[11pt]{article}

\usepackage{amsthm}
\usepackage{amsmath}
\usepackage{amssymb}
\usepackage{latexsym}
\usepackage{graphicx}
\usepackage{subfigure}

\newcommand{\bpr}{{\bf Proof.} \hspace{1 em}}
\newcommand{\epr}{ \\ \hspace*{4.5in} $\Box$ }

\newfam \msbfam
\font \tenmsb=msbm10	\textfont \msbfam=\tenmsb

\newcommand{\Rm}{\mathbb{R}}

\newcommand\bz{{\bf 0}}
\newtheorem{lem}{Lemma}
\newtheorem{thm}{Theorem}

\newtheorem{cor}{Corollary}
\newtheorem{prop}{Proposition}

\newcommand{\tr}{\mbox{Trace}}

\begin{document}

\title{Toward the Universal Rigidity of General Frameworks
\thanks{Alfakih's research is supported by the Natural Sciences and Engineering Research Council of Canada.
Taheri's research is supported in part by DOE Grant DE-SC0002009.
 Ye's research supported in part by
NSF Grant GOALI 0800151 and DOE Grant DE-SC0002009.}
}

\author{Abdo Y. Alfakih
\thanks{Department of Mathematics and Statistics, University of Windsor,
		Windsor, Ontario N9B 3P4, Canada.  Research supported by the Natural Sciences
		and Engineering Research Council of Canada.  E-mail: {\tt
		alfakih@uwindsor.ca}},\quad
		Nicole Taheri \thanks{Institute for Computational and Mathematical
		Engineering, Stanford University, Stanford, CA 94305. Research
		supported in parth by DOE Grant de-sc00020009. E--mail: {\tt
		ntaheri@stanford.edu}},
		\quad and\quad
		Yinyu Ye\thanks{Department of Management Science and Engineering,
		Stanford University, Stanford, CA 94305. E--mail: {\tt
		yinyu-ye@stanford.edu}}}


\date{\today}
\maketitle

\begin{abstract}
Let $(G,P)$ be a bar framework of $n$ vertices in general position in $\Rm^d$,
for $d \leq n-1$, where $G$ is a $(d+1)$-lateration graph. In this paper, we
present a constructive proof that $(G,P)$ admits a positive semidefinite stress
matrix with rank $(n-d-1)$.  We also prove a similar result for a sensor
network, where the graph consists of $m(\ge d+1)$ anchors.
\end{abstract}

\section{Introduction}
Let $V(G)$ and $E(G)$ be, respectively, the vertex set and the edge set of a
simple edge-weighted graph $G$, where each edge $(i,j)$ has a positive weight
$d_{ij}$. The {\em graph realization problem} (GRP) is the problem of
determining whether there exists a realization of $G$ in Euclidean space
$\Rm^d$, for a given dimension $d$.  A (matrix) {\em realization}
$P=[p_1,\ldots,p_n]$ of $G$ in $\Rm^{d}$ is a mapping $P: V(G) \rightarrow
\Rm^{d\times n}$ such that, if vertices $i$ and $j$ of $G$ are adjacent, then
the Euclidean distance between points $p_i\in \Rm^d$ and $p_j\in \Rm^d$ is
equal to the prescribed weight $d_{ij}$ on the edge $(i,j)$. We always assume
that the points $p_1,\ldots,p_n$ affinely span $\Rm^d$. In order words, $P$ is
a realization of $G$ if and only if $p_1,\ldots,p_n$ affinely span $\Rm^d$ and
\[\|p_i-p_j\|=d_{ij} \ \mbox{ for each edge } (i,j) \in E(G).\]
Throughout this paper, $\|x\|$ denotes the $2$-norm of a vector $x$.
Also, we use $\bz$ to denote the matrix of all zeros of the appropriate dimension.
See, e.g.,
\cite{GSS91,H92,A00,A01,JJ03,Connelly,EGW+04,SY05,SY06,A07,Gortler,Zhu}.

The GRP and its variants arise from applications in various areas, such as
molecular conformation, dimensionality reduction, Euclidean ball packing, and
more recently, wireless sensor network localization
\cite{AKW,BY,AGY04,SY05,BTY,SY06,So07}.

Let $P$ be a given realization of graph $G$ with $n$ vertices in $\Rm^d$.  A
realization $P$ together with a graph $G$ is often referred to as a {\em bar
framework} (or {\em framework}), and is denoted by $(G,P)$.  If $P$ is the only
realization of $G$ in $\Rm^d$, up to a rigid motion (e. g., translation or
rotation), then we say that the framework $(G,P)$ is {\em globally rigid}.
However, if $P$ is the only realization of $G$, up to a rigid motion, in all
dimensions, then we say that the framework $(G,P)$ is {\em universally rigid}.

For a given framework $(G,P)$ in $\Rm^d$, define the $(d+1) \times n$ matrix $A$ such that
\begin{equation} \label{defA}  A = \left[ \begin{array}{c} P \\ e^T \end{array} \right] , \end{equation}
where $e$ is the vector of all 1's in $\Rm^n$. Matrices $P$ and $A$ are also respectively called
the {\em position matrix } and the {\em extended position matrix} of the framework $(G,P)$.
The notion of a stress matrix plays a critical role in the characterization of the universal, as well as
the global, rigidity of frameworks.
An $n\times n$ symmetric matrix $S$ is called a {\em stress matrix} of framework $(G,P)$
if and only if
\begin{equation}
	\label{null}
	AS=\bz,
\end{equation}
and
\begin{equation}
	\label{noedge}
	S_{ij}=0,\ \forall (i,j)\not\in E(G),
\end{equation}
where $A$ is the extended position matrix of $(G,P)$.  Note that the highest
possible rank of a stress matrix $S$ is $(n-d-1)$, and the zero matrix is a
trivial stress matrix.

The following theorem characterizes the universal rigidity of generic
frameworks in terms of stress matrices.  A framework $(G,P)$ is said to be {\em
generic}, or in {\em generic position}, if the coordinates of $p_1,\ldots,p_n$
are algebraically independent over the integers, i.e., if there does not exist a non-zero
polynomial $f$ with integer coefficients such that $f(p_1,\ldots,p_n)=0$.
\begin{thm} \label{gur}
Let $(G,P)$ be a framework of $n$ vertices in generic position in $\Rm^d$, $d
\leq n-1$. Then $(G,P)$ is universally rigid if and only if there exists a
stress matrix $S$ of $(G,P)$ such that $S$ is positive semidefinite (PSD) and
the rank of $S$ is $n-d-1$.
\end{thm}
The ``if'' part of this theorem was proved independently in \cite{Connelly99}
and \cite{A07}, while the ``only if" part was proved in \cite{Gortler09}.

One of the major research topics in rigidity is whether a result similar to
Theorem \ref{gur} holds if the assumption of a framework in generic position is
replaced by the weaker assumption of a framework in general position.  We say
that framework $(G,P)$ is in {\em general position} in $\Rm^d$ if no $(d+1)$
points of $p_1,\ldots,p_n$ are affinely dependent. For example, points are in
general position in $\Rm^2$ if no $3$ of them are collinear.  It then easily
follows that if framework $(G,P)$ in $\Rm^d$ is in general position, then every
$(d+1)\times (d+1)$ square sub-matrix of the extended position matrix $A$,
defined in (\ref{defA}), has rank $(d+1)$.  Note that whether or not $n$
rational points are in general position can be checked in time polynomial in
$n$ for any fixed dimension $d$, while the generic position condition is {\em
uncheckable}.  The following theorem, proved in \cite{AY10} recently, shows
that the ``if" part of Theorem \ref{gur} still holds true under the general
position assumption.
\begin{thm} \label{gurs}
Let $(G,P)$ be a framework of $n$ vertices in general position in $\Rm^d$, $d
\leq n-1$. Then $(G,P)$ is universally rigid if there exists a stress matrix
$S$ of $(G,P)$ such that $S$ is positive semidefinite and the rank of $S$
equals $n-d-1$.
\end{thm}

However, it remains an open question whether or not the converse of Theorem
\ref{gurs} holds true.  In this paper, we
settle this question in the {\em affirmative} for frameworks $(G,P)$ in general
position when $G$ is a $(d+1)$-latertation graph.  A graph of $n$ vertices is
called a \emph{$(d+1)$-lateration graph} if there is a permutation $\pi$ of the
vertices, $\pi(1), \pi(2), \ldots, \pi(n)$, such that
\begin{itemize}
	\item the first $(d+1)$ vertices, $\pi(1), \ldots, \pi(d+1)$, form a clique, and
	\item each remaining vertex $\pi(j)$, for $j = (d+2), \ldots n$, is adjacent
		to $(d+1)$ vertices in the set $\{\pi(1), \pi(2), \ldots, \pi(j-1)\}$.
\end{itemize}
Such frameworks were shown to be universally rigid in \cite{So07} and
\cite{Zhu}, where several classes of universally rigid frameworks in general
position were identified.

In particular, we present a {\em constructive} proof that a framework $(G,P)$
of $n$ vertices in general position in $\Rm^d$, where $G$ is a
$(d+1)$-lateration graph, admits a PSD and rank $(n-d-1)$ stress matrix $S$. We
show that such a stress matrix $S$ can be computed in strongly polynomial time,
if the $(d+1)$-lateration ordering is known. We also show that if a graph $G$
contains a $(d+1)$-lateration graph as a spanning subgraph, then the framework
$(G,P)$  in general position also admits a PSD stress matrix of rank $(n-d-1)$.
Finally, a similar result for sensor networks, where the graph consists of
$m(\ge d+1)$ anchors is also given.

\section{The GRP and Semidefinite Programming (SDP)}
\label{pre-stress}
If the graph realization problem is relaxed to the problem of determining
whether a realization of the given edge-weighted graph $G$ exists in some
unspecified Euclidean space, then this relaxed problem can be modeled as a
semidefinite programming problem (SDP).  Furthermore, one can find a stress
matrix with the maximum rank for any given framework by solving a pair of
semidefinite programs (see also \cite{SY05,SY06,A07,BTY}).  In particular, one
can formulate a pair of dual SDPs where $A^TA$ is a solution to the primal
problem, and the stress matrix is a solution to the dual problem. Here, $A$ is
the extended position matrix defined in (\ref{defA}). Next, we present one such
formulation (for other SDP formulations of the same problem, see
\cite{A10,AKW,BY,SY05}).

Let the inner product of two matrices $R$ and $Q$ be defined by $R\cdot
Q=\tr(R^TQ)$.  An SDP for the relaxed graph realization problem attempts to
find a symmetric matrix $Y\in \Rm^{n\times n}$ that solves
\begin{equation}
	\label{edgpy}
	\begin{array}{cl@{\,\,\,}l}
		 \mbox{maximize} & \bz\cdot Y & \\
		 \noalign{\medskip}
		 \mbox{subject to} & (e_{i}-e_j)(e_i-e_j)^T\cdot Y=d_{ij}^2, & \forall\ (i<j,j)\in E(G)\\
		 \noalign{\medskip}
		 & Y\succeq \bz
	\end{array}
\end{equation}
where $e_j\in \Rm^n$ is the vector of all zeros except $1$ at the $j$th
position, and $Y \succeq \bz$ constrains $Y$ to be symmetric PSD.
$A^TA$ and $P^TP$ are both feasible solutions to  Problem (\ref{edgpy})
since
\[(e_{i}-e_j)(e_i-e_j)^T\cdot A^TA=\|a_i-a_j\|^2=\|p_i-p_j\|^2=d_{ij}^2,\ \forall\ (i<j,j)\in E(G).\]
The dual of Problem (\ref{edgpy}) is:
\begin{equation}\label{edgpyd}
\begin{array}{cl}
\mbox{minimize} & \sum_{(i<j,j)\in E(G)}w_{ij}d_{ij}^2 \\
\noalign{\medskip}
\mbox{subject to}& S:=\displaystyle{
\sum_{(i<j,j)\in E(G)}w_{ij}(e_i-e_j)(e_i-e_j)^T\succeq \bz}
\end{array}
\end{equation}
Note that the dual problem is always feasible, since $w_{ij}=0$ for all
$(i,j)\in E(G)$ is a feasible solution. In fact, this solution is also optimal,
since by the weak duality theorem, $0$ is a lower bound on the objective.

From the duality theorem, any optimal solution $S$ of (\ref{edgpyd}) and any
feasible solution $Y$ of (\ref{edgpy}) will satisfy $Y\cdot S=\bz$. This
implies $A^TA\cdot S=ASA^T=0$, or $AS=\bz$.  Moreover, $S_{ij}=0, \; \forall
(i,j)\not\in E(G)$, so that any dual optimal solution is a PSD stress matrix.
We say that the SDPs (\ref{edgpy}) and (\ref{edgpyd}) admit a strictly
complementary solution pair when their respective solutions $(Y,S)$ satisfy
$\mbox{\rm rank}(Y)+\mbox{\rm rank}(S)=n$.

Thus, the question of determining whether there is a non-trivial PSD stress
matrix is equivalent to determining whether there is a non-trivial dual optimal
solution, given that the primal problem is feasible.  In particular, when the
framework is universally rigid in $\Rm^d$, the primal problem (\ref{edgpy}) has
a solution $Y=A^TA$ with rank $(d+1)$. Hence, the SDP Problems (\ref{edgpy})
and (\ref{edgpyd}) admit a strictly complementary solution pair if and only if
there is a dual optimal solution $S$ for (\ref{edgpyd}) with rank $(n-d-1)$.
\begin{prop}
	A universally rigid framework of $n$ vertices in $\Rm^d$, $d \leq n-1$, always
	admits a non-trivial positive semidefinite stress matrix.
\end{prop}
\bpr
	This follows simply from Theorem 6 in \cite{A10}, which states that a
	framework $(G,P)$ in $\Rm^{d}$ admits a non-trivial PSD stress matrix if and
	only if there does not exist a framework $(G,Q)$ in $\Rm^{n-1}$ such that
	$||q_i-q_j||=||p_i-p_j||$ for all $(i,j) \in E(G)$.
\epr

The following result, stated in \cite{SY05}, answers the question
of whether we could find a non-trivial PSD stress matrix, if it exists.
\begin{prop}\label{compute}
	A primal solution $Y$ of (\ref{edgpy}) that has the highest possible rank
	among all primal feasible solutions, together with a dual solution $S$ of
	(\ref{edgpyd}) that has the highest possible rank among all dual optimal
	solutions, can be computed approximately by an SDP interior-point algorithm
	in polynomial time of $n$, $d$, and $\log(1/\epsilon)$ with error $\epsilon$.
\end{prop}
Proposition \ref{compute} also implies that if a universally rigid framework
of $n$ vertices in $\Rm^d$, $d \leq n-1$, admits a rank $(n-d-1)$ and PSD
stress matrix, then such a stress matrix can be computed approximately in
polynomial time. However, we may not be able to compute such a stress matrix
exactly using the SDP algorithm, even when $Y$ is known.


\section{Main Result}
\label{secm}
The following theorem, whose proof is given at the end of this section, is our
main result.

\begin{thm}	\label{main}
Let $(G,P)$ be a framework of $n$ vertices in general position in $\Rm^d$, $d
\leq n-1$, where $G$ is a $(d+1)$-lateration graph. Then $(G,P)$ admits a
positive semidefinite stress matrix with rank $(n-d-1)$. Moreover, such a
stress matrix can be computed {\em exactly} in {\em strongly} polynomial time,
$\mathcal{O}(n^3+nd^3)$ arithmetic operations, if the lateration ordering and
the position matrix $P$ are known.
\end{thm}

An $n \times n$ symmetric matrix $S$ that satisfies condition (\ref{null}),
i.e., $AS = \bz$, is called a \emph{pre-stress matrix}.\footnote{The term
\emph{pre-stress} has been used by Connelly \emph{et al} to mean something
different, see \cite{Connelly92}.} Our constructive proof of Theorem \ref{main}
first generates a PSD pre-stress matrix with rank $(n-d-1)$, then uses this
pre-stress matrix as a basis to generate a PSD stress matrix with rank
$(n-d-1)$. Recall that a stress matrix is a pre-stress matrix which also
satisfies condition (\ref{noedge}), i.e., $S_{ij} = 0, \ \forall (i,j) \not\in
E(G)$.

 The following result follows from basic linear algebra.
\begin{prop}
	\label{stress}
	For any framework in $\Rm^d$, there exists a pre-stress matrix which is
	positive semidefinite and has rank $(n-d-1)$.  Moreover, a universally rigid
	framework in $\Rm^d$ on a complete graph has a rank $(n-d-1)$ positive
	semidefinite stress matrix.
\end{prop}
For example, the projection matrix
\[I-A^T(AA^T)^{-1}A,\]
where $A$ is the extended position matrix,
is a PSD pre-stress matrix with rank $(n-d-1)$. Clearly, the projection matrix can be constructed in
$\mathcal{O}(n^3)$ arithmetic operations.

Under the general position assumption, one can find a matrix $L \in
\Rm^{n\times(n-d-1)}$ of the form
\[ L = \begin{pmatrix} * & * & \cdots & * & * \\
											 \vdots & \vdots & \cdots & \vdots & \vdots \\
											 * & \vdots & \ddots & \vdots & \vdots \\
											 1 & * & \cdots & * & * \\
											 0 & 1 & \cdots & * & * \\
											 \vdots & \vdots & \ddots & \vdots & \vdots \\
											 0 & 0 & \cdots & 1 & * \\
											 0 & 0 & \cdots & 0 & 1 \\
			\end{pmatrix}, \]
that is, for $k=1,...,(n-d-1)$, $L_{ik}=1$ for $i=d+1+k$ and $L_{ik}=0$ for
$i>(d+k+1)$, such that
\[AL=\bz,\]
where $A$ is the extended position matrix.  Clearly, $L$ has rank $(n-d-1)$,
thus $S=LL^T$ is a PSD pre-stress matrix with rank $(n-d-1)$.  Such a matrix
$L$ is called a Gale matrix of framework $(G,P)$ since its columns form a basis
for the nullspace of $A$ \cite{A10}.

For a $(d+1)$-lateration graph $G$ with lateration ordering $1,2,\ldots,n$, and
for a vertex $k \in V(G)$, let
\begin{equation} \label{defNk}
N(k) = \{ i \in V(G): i \leq k-1 \mbox{ and } (i,k) \in E(G) \}.
\end{equation}
Thus, for such a graph, $|N(k)| = d+1$ for each vertex $k=d+2,\ldots,n$.
Furthermore, one can generate the $k$th column of $L$, $L_k$, for
$k=1,...,(n-d-1)$, by solving the system of linear equations
\begin{equation}
	\label{leq1}
	\sum_{i\in N(k)}L_{ik}a_i=-a_{d+k+1},
\end{equation}
where $a_{i}$ is $i$th column of the extended position matrix $A$, and assigning
$L_{ik}=0$ for all $i\not\in N(k)$. The above $d\times d$ linear equation system can be solved
in $\mathcal{O}(d^3)$ operations and there are $n-d-1$ many of them to solve, and the formation of $S$ takes at most $\mathcal{O}(n^3)$
operations. Therefore, we have the following theorem.
\begin{lem}
	\label{Sn}
	The linear system (\ref{leq1}) has a unique solution under the general position
	condition. Moreover, the matrix
		\[S^{n}=LL^T=\sum_{k=1}^{n-d-1}L_kL_k^T\succeq \bz\]
	is a pre-stress matrix with rank $(n-d-1)$, and can be computed in
	$\mathcal{O}(n^3+nd^3)$ arithmetic operations.
\end{lem}

Next, we present an algorithm which uses $S^n$ of Lemma \ref{Sn} as a basis to generate
the desired stress matrix.

\subsection{A Purification Algorithm}

If the pre-stress matrix $S^{n}$, as constructed in Lemma \ref{Sn},
satisfies condition (\ref{noedge}), i.e., $S_{ij} = 0, \ \forall (i,j) \not\in
E(G)$, then it is the desired stress matrix. This is true if the graph is a
$(d+1)$-tree graph, that is, if there is a permutation $\pi$ of the vertices such
that,
\begin{itemize}
	\item the first $(d+1)$ vertices, $\pi(1), \ldots, \pi(d+1)$, form a clique, and
	\item each vertex $\pi(j)$, for $j = (d+2), \ldots, n$, is adjacent to the $(d+1)$ vertices of a
		$(d+1)$-clique in the set $\{\pi(1),\pi(2),\cdots,\pi(j-1)\}$.
\end{itemize}
In this case, any entry in $S^n=LL^T$, for $i<j$ and $(i,j)\not\in E(G)$, is zero.

However, if $S^{n}$ is not a stress matrix, we need to zero out the entries
which should be zero but are not, i.e., the entries $S^n_{ij}\ne 0$, $i<j\
\mbox{and}\ (i,j)\not\in E(G)$.  We do this in reverse order by column; first, we
zero out the entries $S^n_{in}\ne 0$, for $i<n$ and $(i,n)\not\in E(G)$, and then
do the same for columns $(n-1), (n-2), \ldots, (d+3)$.  This ``purification''
process will keep the pre-stress matrix PSD and maintain rank $(n-d-1)$.

If $S^n$ is constructed from $L$ as in the previous section, there is no need for
purification of the last column (or row), since any entry in $LL^T$ for $i<n$
and $(i,n)\not\in E(G)$ is zero. But for general pre-stress matrices,
this may not be the case. Therefore, we first show how to purify the last
column (or row) of a PSD pre-stress matrix with rank $(n-d-1)$.
We construct a vector $s^n\in \Rm^n$ with the elements,
	\[ s^n_{i}=-S^n_{in}\;, \forall (i,n)\not\in E(G) \quad \mbox{and} \quad s^n_{n}=1, \]
and solve the following system of linear equations for the remaining entries in $s^n$,
\begin{equation}
	\label{leq2}
	\sum_{i\in N(n)}s^n_ia_i=-\sum_{(i,n)\not\in E(G)}s^n_ia_i.
\end{equation}
The right-hand-side of the equation can be formed in at most $\mathcal{O}(nd)$ operations, and
the $d\times d$ linear system can be solved in $\mathcal{O}(d^3)$ operations. Thus, $s^n$ can be
computed in at most $\mathcal{O}(nd+d^3)$
operations.

The linear system (\ref{leq2}) has a unique solution under the
general position condition, and by construction, $As^n=\bz$.
\begin{lem}
	\label{stepn}
	Let $S^{n-1}=S^n+s^n(s^n)^T$. Then
	\begin{itemize}
		\item $AS^{n-1}=\bz$.
		\item $S^{n-1}\succeq \bz$ and the rank of $S^{n-1}$ remains $(n-d-1)$.
		\item $S^{n-1}_{in}=0$ for all $i<n,\ (i,n)\not\in E(G)$.
	\end{itemize}
\end{lem}
\bpr
The first statement holds, since
	\[AS^{n-1}=AS^n+As^n(s^n)^T=As^n(s^n)^T = \bz,\]
where the last step follows from the construction of $s^n$, so that $As^n=\bz$.

The second statement follows from $S^{n-1} = S^n + s^n (s^n)^T \succeq
S^n\succeq \bz$. Thus, $\mbox{rank}(S^{n-1})\ge \mbox{rank}(S^{n})=(n-d-1)$, but
$AS^{n-1}=\bz$ implies that the rank of $S^{n-1}$ is bounded above by $(n-d-1)$.

The third statement is also true by construction. In the last column (or row)
of $s^n(s^n)^T$, the $i$th entry, where $i\ne n$ and $(i,n)\not\in E(G)$, is
precisely $-S^n_{in}$, i.e.,
	\[ \left( s^n (s^n)^T \right)_{in} = s^n_i s^n_n = s^n_i =  -S^n_{in}, \]
so that it is canceled out in the last column (or row) of matrix
$S^{n-1}=S^n+s^n(s^n)^T$.
\epr
Note that update $S^{n-1}=S^n+s^n(s^n)^T$ uses $\mathcal{O}(n^2)$ arithmetic operations.

We continue this purification process for $(n-1),\ldots,k, \ldots, (d+3)$.
Before the $k$th purification step, we have $S^{k}\succeq \bz$,
$AS^k=\bz$, rank$(S^k) = (n-d-1)$, and
	\[S^{k}_{ij}=0,\ \forall j>k, i<j\ \mbox{and}\ (i,j) \not\in E(G)\]
We then construct a vector $s^k\in \Rm^n$ with the elements,
	\[s^k_{i}=-S^k_{ik}, \ \forall (i,k) \not\in E(G),\ s^k_{k}=1, \quad \mbox{and} \quad s^k_i=0\ \forall i>k,\]
and solve the system of linear equations for the remaining entries in $s^k$:
\begin{equation}
	\label{leq3}
	\sum_{(i,k)\in E(G)}s^k_ia_i=-\sum_{(i,k)\not\in E(G)}s^k_ia_i.
\end{equation}
Again, solving this linear system takes at most $\mathcal{O}(nd+d^3)$ operations,
and by construction, we have $As^k=\bz$.

Similarly, the following lemma shows results analogous to those in Lemma
\ref{stepn}, for the remaining columns.
\begin{lem}
	\label{stepk}
	Let $S^{k-1}=S^k+s^k(s^k)^T$. Then
	\begin{itemize}
		\item $AS^{k-1}=\bz$.
		\item $S^{k-1}\succeq \bz$ and the rank of $S^{k-1}$ remains $(n-d-1)$.
		\item $S^{k-1}_{ij}=0$ for all $j\ge k$ and $i< j,\ (i,j)\not\in E(G)$.
	\end{itemize}
\end{lem}
\bpr
The proof of the first two statements is identical to that in Lemma \ref{stepn}.

The third statement is again true by construction. Note that in the $k$th
column (or row) of $s^k(s^k)^T$, the $i$th entry, $i>k$ and $(i,k)\not\in E(G)$, is
precisely $-S^k_{ik}$, i.e.,
	\[ \left( s^k (s^k)^T \right)_{ik} = s^k_i s^k_k = s^k_i =  -S^k_{ik}, \]
so that it is canceled out in the $k$th column (or row) of matrix
$S^{k-1}=S^k+s^k(s^k)^T$. Furthermore, for $j = (k+1), \ldots, n$, the $j$th
column (or row) of $s^k(s^k)^T$ has all zero entries, which means the entries
in $j$th column (or row) of $S^{k-1}$ remain unchanged from $S^k$.
\epr

Now we are ready to prove our main result.

\noindent {\bf Proof of Theorem \ref{main}.}

    Assume that the $(d+1)$-lateration graph has the lateration ordering $1,2,\ldots,n$.
	The matrix $S^{d+2}$, constructed via the process described in Lemmas
	\ref{stepn} and \ref{stepk}, will be a PSD stress matrix with rank $(n-d-1)$,
	for the $(d+1)$-lateration graph, since after step $k=(d+3)$, we will have a
	matrix $S^{d+2}$ that satisfies,
		\[ AS^{d+2}=\bz \quad \mbox{and} \quad S^{d+2}_{ij}=0, \; \forall (i,j) \not\in E(G) \]
	Note that the first $(d+2)$ vertices form a clique in $G$, and the
	principal $(d+2)\times(d+2)$ submatrix has no zero entries. This stress
	matrix is unique and always exists since the graph is a $(d+1)$-lateration
	graph, and thus there is always a unique solution to the linear equation
	(\ref{leq3}).  Furthermore, by Lemma \ref{stepk}, $S^{d+2}\succeq \bz$ and
	the rank of $S^{d+2}$ remains $(n-d-1)$.

	There are $(n-d-2)$ purification steps, where each step computes a rank-one
	matrix $s^k(s^k)^T$ and forms a new pre-stress matrix $S^k+s^k(s^k)^T$, taking
	at most $\mathcal{O}(n^2+d^3)$ arithmetic operations. Thus, the computation of the
      max-rank PSD stress matrix uses at most $\mathcal{O}(n^3+nd^3)$ operations.
\epr

We also have the following corollary:
\begin{cor}
	\label{main2}
	Any universally rigid framework $(G,P)$ in general position admits a
	positive semidefinite stress matrix with rank $(n-d-1)$, if $G$ contains a
	$(d+1)$-lateration graph as a spanning subgraph. Moreover, such a stress
	matrix can be computed {\em exactly} in {\em strongly} polynomial time,
	$\mathcal{O}(n^3+nd^3)$ arithmetic operations, if the lateration ordering and
	the position matrix $P$ are known. Otherwise, such a stress matrix, together
	with the position matrix $P$, can be computed approximately by an SDP
	interior-point algorithm in time polynomial in $n$, $d$, and
	$\log(1/\epsilon)$, with error $\epsilon$.
\end{cor}

This secondary result holds because we can ignore all edges outside of the
$(d+1)$-lateration spanning subgraph to prove the existence of a PSD stress
matrix with rank $(n-d-1)$. Since finding a $(d+1)$-lateration spanning
subgraph requires at least $\mathcal{O}(n^{d+2})$ operations, we cannot
actually construct such a stress matrix exactly in $\mathcal{O}(n^3+nd^3)$
operations, if either the lateration ordering or the position matrix $P$ is
unknown. However, Proposition \ref{compute} implies that such a rank $(n-d-1)$
and PSD stress matrix, together with the position matrix $P$, can be computed
approximately in polynomial time, although not strongly polynomial.

\section{Strong Localizability of $(d+1)$-Lateration Graph with Anchors}
In this section we study the stress matrix of a \emph{sensor network}, or
graph localization with anchors.  A sensor network consists of $m(\ge d+1)$ anchor points whose
positions, $\bar{p}_1,\ldots,\bar{p}_{m}\in \Rm^d$, are known, and $n$ sensor
points whose locations, $x_1,\ldots,x_n\in \Rm^d$, are yet to be determined.
We are given the Euclidean distance values $\bar{d}_{kj}$ between
$\bar{p}_k$ and $x_j$ for some $(k,j)$, and $d_{ij}$ between $x_i$ and $x_j$ for
some $i<j$.  Specifically, let
	\[ N_a=\{(k,j):\bar{d}_{kj} \mbox{ is specified}\} \quad \mbox{and} \quad
	   N_x=\{(i,j):i<j,\, d_{ij} \mbox{ is specified}\}. \]
The problem is to find a realization of $x_1,\ldots,x_n\in \Rm^d$ such that
\begin{equation}
	\label{edgp}
	\begin{array}{ll}
		\|\bar{p}_k-x_j\|^2 = \bar{d}_{kj}^2 & \forall\ (k,j)\in N_a\\
		\noalign{\medskip}
		\|x_i-x_j\|^2 = d_{ij}^2 & \forall\ (i,j)\in N_x.
	\end{array}
\end{equation}

The semidefinite programming relaxation model for (\ref{edgp}) attempts to find a
$(d+n)\times (d+n)$ symmetric matrix
\begin{equation}
	\label{Zrep}
	Z = \begin{pmatrix} I_d & X \\ X^T & Y \end{pmatrix} \succeq \bz
\end{equation}
that solves the SDP
\begin{equation}\label{edgpz}
\begin{array}{cl@{\,\,\,}l}
   \mbox{maximize} & \bz\cdot Z & \\
   \noalign{\medskip}
   \mbox{subject to} & Z_{1:d,1:d}=I_d & \\
   \noalign{\medskip}
   & (\bz;e_{i}-e_j)(\bz;e_i-e_j)^T\cdot Z=d_{ij}^2 & \forall\ (i,j)\in N_x\\
   \noalign{\medskip}
   & (-\bar{p}_k;e_{j}) (-\bar{p}_k;e_{j})^T\cdot Z = \bar{d}_{kj}^2 & \forall\ (k,j)\in N_a\\
   \noalign{\medskip}
   & Z\succeq \bz,
\end{array}
\end{equation}
where $(-\bar{p}_k;e_j)\in \Rm^{d+n}$ is the vector of $-\bar{p}_k$ vertically
concatenated with $e_j$.  $Z_{1:d,1:d}$ is the $d\times d$ top-left principal submatrix
of $Z$ and $I_d$ is the $d$--dimensional identity matrix.  $Z_{1:d,1:d}=I_d$
can be represented as $d(d+1)/2$ linear equality constraints.

The dual of the SDP relaxation model is given by:
\begin{equation}\label{edgpzd}
	\begin{array}{cl}
		\mbox{minimize} & \displaystyle{I_d\cdot V+\sum_{(i,j)\in N_x}w_{ij}d_{ij}^2+\sum_{(k,j)\in N_a}\bar{w}_{kj}\bar{d}_{kj}^2} \\
		\noalign{\medskip}
		\mbox{subject to}& S:=\displaystyle{\left(\begin{array}{cc}
		V & \bz\\
		\bz & \bz
		\end{array}\right)+\sum_{(i,j)\in N_x}w_{ij}(\bz;e_i-e_j)(\bz;e_i-e_j)^T} \\
		\noalign{\medskip}
		 & \quad \displaystyle{+\sum_{(k,j)\in N_a}\bar{w}_{kj}(-\bar{p}_k;e_{j})(-\bar{p}_k;e_{j})^T\succeq \bz}.
	\end{array}
\end{equation}
Note that the dual is always feasible, since the symmetric matrix with
$V=\bz\in \Rm^{d\times d}$, variables $w_{ij}=0$ for all $(i,j)\in N_x$ and
$\bar{w}_{kj}=0$ for all $(k,j)\in N_a$, is feasible for the dual. Also, each
column $j=(d+1),\ldots,(d+n)$ of the dual matrix $S$ has the structure:
\begin{equation}\label{diags}
\begin{array}{rcl}
S_{1:d,j}&=&-\sum_{(k,j)\in N_a}\bar{w}_{kj}\bar{p}_k,\\
S_{ij}&=&-w_{ij},\  (i,j)\in N_x,\\
S_{ij}&=&0,\ (i,j)\not\in N_x,\\
S_{jj}&=&\sum_{(i,j)\in N_x}w_{ij}+\sum_{(k,j)\in N_a}\bar{w}_{kj}.
\end{array}
\end{equation}
Let $w_{ij}$ and $\bar{w}_{kj}$ be called stress variables, and $S$ be called a
variable stress matrix for the sensor network localization problem.

It is shown in \cite{SY05} that both SDPs (\ref{edgpz}) and (\ref{edgpzd}) are
feasible and solvable when there is at least one anchor point, the graph is
connected, and there is no duality gap between the two SDPs. Let $P$ be a
position matrix of the $n$ sensors satisfying constraints in (\ref{edgp}).
Then, the sensor network is said to be \emph{uniquely localizable} if
\begin{equation}\label{Z-P}
Z= \begin{pmatrix} I_d & P\\ P^T & P^TP \end{pmatrix} \succeq \bz
\end{equation}
is the only matrix solution to the primal SDP  (\ref{edgpz}); this is similar to the
concept of universal rigidity.  The network is said to be \emph{strongly localizable}
if there is an optimal dual stress matrix $S$ such that
\begin{itemize}
	\item $ZS=\bz$,
	\item $S\succeq \bz$ and rank$(S) = n$.
\end{itemize}
It has been shown in \cite{SY05} that strong localizability implies unique localizability.

The standard graph realization problem is equivalent to the sensor network
localization problem without anchors; thus, the two problems are different, but
closely related. For example, unlike in the SDP (\ref{edgpy}), $Z$ constructed
from $A$, where $A$ is the extended position matrix, is no longer feasible for
(\ref{edgpz}), although $Z$ constructed from a position matrix $P$ in
(\ref{Z-P}) is feasible. Hence, the stresses of the dual on the anchors may not
need to be balanced. As another example, consider a sensor network of two
anchors and one sensor in $\Rm^2$, where the distances from the sensor to the
two anchors are known. The network is not uniquely localizable, but it is
universally rigid in graph realization, since the three points form a clique.

However, if the sensor network has at least $(d+1)$ anchors in general
position, and the graph realization problem has a $(d+1)$-point clique also in
general position, then unique localizability is equivalent to universal
rigidity, and strong localizability is equivalent to a framework on $(n+d+1)$
points having a PSD stress matrix with rank $n$ (see \cite{So07,Zhu}). The
latter implies that the SDP pair (\ref{edgpz}) and (\ref{edgpzd}) admits a
strictly complementary solution pair.

\begin{thm}\label{main3}
Take a graph $G$ of $m(\ge d+1)$ anchor points and $n$ sensor points with edges
given in $N_x$ and $N_a$, and let $G$ be a $(d+1)$-lateration graph with
$(\bar{P},P)$ in general positions. Then the sensor network is strongly
localizable, and a rank $n$ optimal dual stress matrix can be computed exactly
in strongly polynomial time, $\mathcal{O}(n^3+nd^3)$ arithmetic operations, if
the lateration ordering and the sensor position matrix $P$ are known.
\end{thm}

\bpr

We need to show that,  in $\mathcal{O}(n^3+nd^3)$ arithmetic operations,  one
can compute a symmetric matrix ${S \in \Rm^{(d+n) \times (d+n)}}$ which
satisfies $ZS =\bz,\ S \succeq \bz$, $\mbox{rank}(S) = n$, and meets the
structure condition (\ref{diags}). The proof is more complicated than that of
Theorem \ref{main}, since anchor positions appear explicitly in the dual stress
matrix.

For simplicity and without loss of generality, we assume there are exactly
$(d+1)$ anchors which are the first $(d+1)$ points in the lateration ordering;
all other points are sensors and ordered $1,\cdots,n$.
Given a position matrix $P$, the primal feasible solution matrix $Z$ in
(\ref{Z-P}) can be written as
\[
Z = \begin{pmatrix} I_d & P \\ P^T & P^T P \end{pmatrix} =
\begin{pmatrix} I_d \\ P^T \end{pmatrix}
\begin{pmatrix} I_d & P \end{pmatrix}
\]
so that the matrix
$L = [-P ;\ \ I_n] \in \Rm^{(d+n) \times n}$ is in the nullspace of $Z$, or
matrix $[I_d\ \ P]$. Moreover, the matrix
	\[
		S^n = L L^T = \begin{pmatrix} P P^T & -P \\  -P^T & I_n \end{pmatrix}\succeq \bz,
	\]
will also be in the nullspace of $Z$, where $\mbox{rank}(S^n) =\mbox{rank}(L) =
n$. One may call $S^n$ a pre-stress matrix for the sensor localization problem.
$S^n$ may not be a true optimal stress matrix since it may not meet the
structure condition (\ref{diags}).

We now modify the elements of $S^n$ in $\mathcal{O}(n^3)$ operations, so that
these conditions are maintained, and the resulting matrix becomes feasible to
the dual problem (\ref{edgpzd}). Note that the elements in the bottom right $n
\times n$ submatrix of any dual feasible $S$, which corresponds to the sensor
to sensor edges, will be
\[
	S_{ij} = \begin{cases} - w_{ij}, & (i,j) \in N_x \\
						  0, & (i,j) \not\in N_x \\
							\sum_{(i,j) \in N_x} w_{ij} + \sum_{(k,j) \in N_a} \bar{w}_{kj}, & i=j.
           \end{cases}
\]
Thus, for each sensor $j$, the $j$th diagonal element is the negative sum of
the edge weights of $w_{ij}$ and $\bar{w}_{kj}$.  Since there are no
constraints on $V$ in the dual (\ref{edgpzd}), any principal $d \times d$
submatrix is feasible as long as it remains positive semidefinite, and the
submatrix corresponding to the sensor to anchor edges is feasible if, for $i
\leq d, j > d$,
\[
	S_{i,j} = - \sum_{(k,j) \in N_a} \bar{w}_{kj} (\bar{p}_k)_i.
\]

Similar to the constructed proof of Theorem \ref{main}, while maintaining its
rank $n$ and keeping it PSD, we modify each column of the pre-stress matrix
$S^n$, starting with the last, $(d+n)$, and continuing to column $(d+1)$, to
make it a true stress matrix, optimal for the dual.  Each column of the matrix
will be modified in the same way.

We modify each column $(d+\ell)$ of the original matrix by constructing a
vector $s^\ell$ so that the the $(d+\ell)$th column of the matrix $S^{\ell -1}
= S^\ell + s^\ell (s^\ell)^T$ is dual feasible.  Moreover, $S^{\ell-1}$ will be
a positive semidefinite matrix with rank $n$, and $ZS^{\ell -1} = 0$.  Again,
when we modify column $d+\ell$, the $(d+j)$th column (or row) of the modified
pre-stress matrix is unchanged for all $j> \ell$.

More precisely, when modifying the $\ell$th column, we construct a vector
$s^{\ell}\in\Rm^{d+n}$ such that $s^{\ell}_{d+\ell}=1$, $s^{\ell}_{i}=0$ for
$i>(d+\ell)$, and the first $(d+\ell-1)$ entries are
\begin{equation}
	\label{sj-eqn}
	s^{\ell}_{1:(d+\ell-1)} := \begin{pmatrix} - \sum_{(k,\ell) \in N_a} \bar{w}_{k\ell} \bar{p}_k \\
                          -\sum_{(i<\ell,\ell) \in N_x} w_{i\ell}e_i \end{pmatrix} - S^{\ell}_{1:(d+\ell-1),(d+\ell)},
\end{equation}
where the $(d+1)$ stress variables $\bar{w}_{k\ell}$ and $w_{i\ell}$ are yet to
be determined, and $e_i\in \Rm^n$ is the vector of all zeros except $1$ at the
$i$th position.

For the updated matrix $S^{\ell-1}:=S^{\ell}+s^{\ell}(s^{\ell})^T$, adding
$s^{\ell}(s^{(\ell)} )^T$ to $S^{\ell}$ will not affect any column (or row) to
the right (or below) of column (or row) $d+\ell$. In particular, the
$(d+\ell)$th column of $S^{\ell-1}$ becomes
\[
S^{\ell-1}_{d+\ell} =
	\begin{pmatrix} S^{\ell-1}_{1:(d+\ell-1),(d+\ell)} \\ \\
		S^{\ell-1}_{ (d+\ell),(d+\ell) } \\ \\
 		S^{\ell-1}_{(d+\ell+1):(d+n),(d+\ell)}
	\end{pmatrix}
=
	\begin{pmatrix}
	\begin{pmatrix} - \sum_{(k,\ell) \in N_a} \bar{w}_{k\ell} \bar{p}_k \\
                -\sum_{(i<\ell,\ell) \in N_x} w_{i\ell}e_i \end{pmatrix} \\ \\
	1+S^{\ell}_{(d+\ell),(d+\ell)} \\ \\
	S^{\ell}_{(d+\ell+1):(d+n),(d+\ell)}
	\end{pmatrix}
\]
In the case where $\ell = n$, there is no last term and this becomes
\[
S^{n-1}_{d+n} =
	\begin{pmatrix}
	\begin{pmatrix} - \sum_{(k,n) \in N_a} \bar{w}_{kn} \bar{p}_k \\
                -\sum_{(i<n,n) \in N_x} w_{in}e_i \end{pmatrix} \\ \\
								1+S^{n}_{(d+n),(d+n)}
	\end{pmatrix}
\]

By construction, column $(d+\ell)$ of $S^{\ell -1}$ almost meets the the
structure conditions of (\ref{diags}).  To ensure $S^{\ell-1}$ is orthogonal to
$Z$, or $[I_d\ \ P]$, we determine the $(d+1)$ stress variables
$\bar{w}_{k\ell}$ and $w_{i\ell}$ in $s^{\ell}$ such that $[I_d \ \ P] s^{\ell}
= \bz$, or equivalently,
\begin{equation}
- \sum_{(k,\ell) \in N_a} \bar{w}_{k\ell}\bar{p}_k - \sum_{(i<\ell,\ell) \in N_x} w_{i\ell}p_i +p_{\ell}
= S^{\ell}_{1:d,(d+\ell)}+\sum_{i=1}^{\ell-1}S^{\ell}_{d+i,(d+\ell)}p_i.
\label{null-j}
\end{equation}
Finally, to meet the diagonal entry value condition of  (\ref{diags}), that is,
the sum of the total edge stresses of sensor $\ell$ equals to the value of its
diagonal element, we add
\begin{equation}
\sum_{(k,\ell) \in N_a}\bar{w}_{k\ell} + \sum_{(i<\ell,\ell) \in N_x}w_{i\ell} = 1 + S^{\ell}_{(d+\ell),(d+\ell)}+\sum_{i=\ell+1}^{n}S^{\ell}_{d+i,(d+\ell)}.
\label{stress-j}
\end{equation}
Equations (\ref{null-j}) and (\ref{stress-j}) are exactly $(d+1)$ linearly
independent equations on the $(d+1)$ stress variables, and thus there is always
a unique solution.

The modified $\ell$th column of $S^{\ell-1}$ is now feasible for a dual
solution stress matrix, and $S^{\ell-1}$ satisfies
\begin{itemize}
	\item $Z S^{\ell-1} = ZS^{\ell} + Zs^{\ell} (s^{\ell})^T = \bz$,
	\item $S^{\ell-1} = S^{\ell} + s^{\ell} (s^{\ell})^T \succeq S^{\ell} \succeq \bz$,
		and hence rank$(S^{\ell-1}) = n$.
\end{itemize}

Repeating this process on each column, for columns $(d+n), \ldots, (d+1)$,
will result in an optimal dual solution matrix $S^0$ that satisfies
\begin{itemize}
	\item $Z S^0 = \bz$,
	\item $S^0 \succeq 0$, and rank$(S^0) = n$.
\end{itemize}

That is, $S^0$ is now a true optimal dual stress matrix with rank $n$. Note
that there are a total of $n$ modification steps, and each modification
step takes at most $\mathcal{O}(n^2+d^3)$ arithmetic operations.
Therefore, we have proved Theorem \ref{main3}.
\epr

Similar to the secondary result for the standard graph realization problem,
we have the following corollary.
\begin{cor}
	\label{main4}
	Take a graph $G$ of $m(\ge d+1)$ anchor points and $n$ sensor points with
	edges given in $N_x$ and $N_a$, and let $G$ contain a $(d+1)$-lateration
	spanning subgraph with $(\bar{P},P)$ in general positions. Then the sensor
	localization problem on $G$ is strongly localizable. Moreover, a rank $n$
	optimal dual stress matrix can be computed {\em exactly} in {\em strongly}
	polynomial time, $\mathcal{O}(n^3+nd^3)$ arithmetic operations, if the
	lateration ordering and the position matrix $P$ are known. Otherwise, such a
	rank $n$ stress matrix, together with the position matrix $P$, can be
	computed approximately by an SDP interior-point algorithm in time polynomial in
	$n$, $d$, and $\log(1/\epsilon)$, with error $\epsilon$.
\end{cor}

The argument for Corollary \ref{main4} is analogous to that of Corollary \ref{main2}.

\section{Examples}
Consider a $3$-lateration framework in dimension 2 on $n=7$ nodes, with ordering
$1,2,\ldots, 7,$ and position matrix
\[P=\left(\begin{array}{rrrrrrr}
                -1 & 1  & 0 & 2 & 1 & -1 & -2\\
                 1 & 1  &0.5& 0 &-1 & -1 &  0
\end{array}\right)\in \Rm^{2\times 7}.
\]

\noindent{\bf Example 1} Let
\[N(4)=\{1,2,3\},\ N(5)=\{1,3,4\},\ N(6)=\{1,2,4\},\ N(7)=\{3,4,5\}.\]
For this example,
\[L = \begin{pmatrix}
    1.5000  &  5.0000 &  -2.0000   &    0\\
   -0.5000  &      0  & 3.0000     &    0\\
   -2.0000  & -8.0000 &        0   &-1.6000\\
    1.0000  &  2.0000 & -2.0000    &1.4000\\
         0  &  1.0000 &       0    &-0.8000\\
         0  &      0  & 1.0000     &    0\\
         0  &       0 &       0    &1.0000
\end{pmatrix}
\]
and we have the pre-stress matrix $S^7=LL^T,$
\[S^7 = \begin{pmatrix}
   31.2500 &  -6.7500 & -43.0000 &  15.5000  &  5.0000 &  -2.0000 &        0\\
   -6.7500 &   9.2500 &   1.0000 &  -6.5000  &       0 &   3.0000 &        0\\
  -43.0000 &   1.0000 &  70.5600 & -20.2400  & -6.7200 &        0 &  -1.6000\\
   15.5000 &  -6.5000 & -20.2400 &  10.9600  &  0.8800 &  -2.0000 &   1.4000\\
    5.0000 &        0 &  -6.7200 &   0.8800  &  1.6400 &        0 &  -0.8000\\
   -2.0000 &   3.0000 &        0 &  -2.0000  &       0 &   1.0000 &        0\\
         0 &        0 &  -1.6000 &   1.4000  & -0.8000 &        0 &   1.0000
\end{pmatrix}.
\]
Note that $S^7$ is already a stress matrix that meets condition
(\ref{noedge}), so that no ``purification'' algorithm is needed. This example
is not interesting, since the graph is actually a $3$-tree graph.

\vskip 0.2in
\noindent{\bf Example 2} Let
\[N(4)=\{1,2,3\},\ N(5)=\{1,3,4\},\ N(6)=\{2,4,5\},\ N(7)=\{1,3,6\}.\]
In this example,
\[L = \begin{pmatrix}
    1.5000  &  5.0000 &        0   &-1.2500\\
   -0.5000  &      0  & -1.0000    &    0  \\
   -2.0000  & -8.0000 &        0   & 1.0000\\
    1.0000  &  2.0000 &  2.0000    & 0     \\
         0  &  1.0000 & -2.0000    & 0     \\
         0  &      0  & 1.0000     &-0.7500  \\
         0  &       0 &       0    &1.0000
\end{pmatrix}
\]
and we the have the pre-stress matrix $S^7=LL^T,$
\[S^7 = \begin{pmatrix}
   28.8125 &  -0.7500 & -44.2500 &  11.5000  &  5.0000  &  0.9375  & -1.2500\\
   -0.7500 &   1.2500 &   1.0000 &  -2.5000  &  2.0000  & -1.0000  &       0\\
  -44.2500 &   1.0000 &  69.0000 & -18.0000  & -8.0000  & -0.7500  &  1.0000\\
   11.5000 &  -2.5000 & -18.0000 &   9.0000  & -2.0000  &  2.0000  &       0\\
    5.0000 &   2.0000 &  -8.0000 &  -2.0000  &  5.0000  & -2.0000  &       0\\
    0.9375 &  -1.0000 &  -0.7500 &   2.0000  & -2.0000  &  1.5625  & -0.7500\\
   -1.2500 &        0 &   1.0000 &        0  &       0  & -0.7500  &  1.0000
\end{pmatrix}.
\]
While the last column (or row) of $S^7$ meets condition (\ref{noedge}), the
rest does not satisfy (\ref{noedge}). We start the purification process from
$k=6$, where $S^6=S^7$. The column vector $s^6$ is generated by first assigning
\[s^6_{1}=-S^6_{1,6}=-0.9375,\ s^6_3=-S^6_{3,6}=0.75,\ s^6_6=1,\ s^6_7=0,\]
and then solving for $(s^6_2,s^6_4,s^6_5)$ from the linear system (\ref{leq3}) to get
\[s^6=\begin{pmatrix}
   -0.9375\\
   -0.0625\\
    0.7500\\
    0.8750\\
   -1.6250\\
    1.0000\\
         0
\end{pmatrix}.
\]
and
$S^5=S^6+s^6(s^6)^T,$
\[
S^5 = \begin{pmatrix}
   29.6914  & -0.6914 & -44.9531  & 10.6797  &  6.5234  &       0  & -1.2500\\
   -0.6914  &  1.2539 &   0.9531  & -2.5547  &  2.1016  & -1.0625  &       0\\
  -44.9531  &  0.9531 &  69.5625  &-17.3438  & -9.2188  &       0  &  1.0000\\
   10.6797  & -2.5547 & -17.3438  &  9.7656  & -3.4219  &  2.8750  &       0\\
    6.5234  &  2.1016 &  -9.2188  & -3.4219  &  7.6406  & -3.6250  &       0\\
         0  & -1.0625 &        0  &  2.8750  & -3.6250  &  2.5625  & -0.7500\\
   -1.2500  &       0 &   1.0000  &       0  &       0  & -0.7500  &  1.0000
\end{pmatrix}.
\]
Next the column vector $s^5$ is generated by first assigning
\[s^5_{2}=-S^5_{2,5}=-2.1016,\ s^5_5=1,\ s^5_6=s^5_7=0,\]
and solving for $(s^5_1,s^5_3,s^5_4)$ from linear system (\ref{leq3}),
\[s^5=\begin{pmatrix}
   11.3047\\
   -2.1016\\
  -16.4063\\
    6.2031\\
    1.0000\\
         0\\
         0\\
\end{pmatrix}
\]
and
$S^4=S^5+s^5(s^5)^T$,
\[
S^4 = \begin{pmatrix}
  157.4874 & -24.4489 &-230.4207 &  80.8041 &  17.8281  &       0  & -1.2500\\
  -24.4489 &   5.6705 &  35.4319 & -15.5909 &        0  & -1.0625  &       0\\
 -230.4207 &  35.4319 & 338.7275 &-119.1138 & -25.6250  &       0  &  1.0000\\
   80.8041 & -15.5909 &-119.1138 &  48.2444 &   2.7813  &  2.8750  &       0\\
   17.8281 &        0 & -25.6250 &   2.7813 &   8.6406  & -3.6250  &       0\\
         0 &  -1.0625 &        0 &   2.8750 &  -3.6250  &  2.5625  & -0.7500\\
   -1.2500 &        0 &   1.0000 &        0 &        0  & -0.7500  &  1.0000
\end{pmatrix}.
\]
One can see that $S^4$ is now a desired stress matrix for Example 2.


\section*{Acknowledgments}
The authors would like to thank two anonymous referees for their helpful comments
and quick response.

\end{document}